\documentclass[12pt]{amsart}

\usepackage{amssymb}

\newtheorem{theorem}{Theorem}[section]
\newtheorem{proposition}[theorem]{Proposition}
\newtheorem{definition}[theorem]{Definition}
\newtheorem{remark}[theorem]{Remark}

\numberwithin{equation}{section}

\begin{document}

\baselineskip=15.5pt

\title{Connections on a parabolic principal bundle, II}

\author{Indranil Biswas}

\address{School of Mathematics, Tata Institute of Fundamental
Research, Homi Bhabha Road, Bombay 400005, India}

\email{indranil@math.tifr.res.in}



\date{}

\begin{abstract}

In \cite{Bi2} we defined connections on a parabolic principal
bundle. While connections on usual principal bundles are
defined as splittings of the Atiyah exact sequence, it was noted in
\cite{Bi2} that the Atiyah exact sequence does not generalize to
the parabolic principal bundles.
Here we show that a twisted version
of the Atiyah exact sequence generalize to the context of
parabolic principal bundles. For usual principal bundles, giving a
splitting of this twisted Atiyah exact sequence is equivalent
to giving a splitting of the Atiyah exact sequence. Connections on
a parabolic principal bundle can be defined using the
generalization of the twisted Atiyah exact sequence.

\end{abstract}

\maketitle

\section{Introduction}

Generalizing the notion of a parabolic vector bundle, in \cite{BBN1}
and \cite{BBN2} the notion of a parabolic principal bundle was
introduced. Let $G$ be a connected complex linear algebraic group.
A parabolic $G$--bundle over a complex
smooth projective variety $X$
is a smooth variety $E_G$ equipped with an action of $G$ 
as well as a projection to $X$
such that $E_G$ is a principal bundle over the complement
of a simple normal crossing divisor in $X$. However the action of
$G$ over the divisor is allowed to have finite isotropies.

In \cite{Bi2}, connections on a parabolic principal
bundle were defined. Before we describe connections on a parabolic
principal bundle, we will first
briefly recall the definition of a connection
on an usual principal $G$--bundle. Let $F_G$ be a holomorphic
principal $G$--bundle over a complex manifold $Y$, and let
\begin{equation}\label{eq-1}
0\, \longrightarrow\, \text{ad}(F_G)\, \longrightarrow\,
\text{At}(F_G) \longrightarrow\, TY \longrightarrow\,0
\end{equation}
be the corresponding Atiyah exact sequence over $Y$.
A holomorphic (respectively,
complex) connection on $F_G$ is a holomorphic (respectively,
$C^\infty$) splitting of the above Atiyah exact sequence.
Giving a holomorphic (respectively, $C^\infty$) splitting
of \eqref{eq-1} is equivalent to giving a holomorphic (respectively,
$C^\infty$) one--form $\omega$ on $F_G$ with values in the Lie
algebra $\mathfrak g$ of $G$ and satisfying the following two
conditions:
\begin{itemize}
\item the restriction of $\omega$ to each fiber of the projection
$F_G\, \longrightarrow\, Y$ coincides with the Maurer--Cartan form,
and
\item the form $\omega$ is $G$-equivariant for the adjoint action
of $G$ on $\mathfrak g$.
\end{itemize}

Let $E_G$ be a parabolic $G$--bundle over $X$. Then there exists
a Galois covering $f\,:\, Y\, \longrightarrow\, X$, where $Y$
is a smooth complex variety, together with
a holomorphic principal $G$--bundle
$F_G$ over $Y$ equipped with a lift of the action of the Galois
group $\Gamma\, :=\, \text{Gal}(f)$ on $Y$, such that $E_G\,
=\, F_G/\Gamma$. It should be mentioned that there are many
coverings of $X$ satisfying the above conditions.
In \cite{Bi2} we noted that there is no Atiyah
exact sequence for a general parabolic $G$--bundle. This means that
for a choice of a covering $Y$ of the above type, the exact
sequence of vector bundles over $X$ given by the Atiyah exact
sequence for $F_G$ on $Y$ depends on the choice of covering.
Since the Atiyah exact sequence is not available, we used the
above description of a connection as a $\mathfrak
g$--valued $1$--form on the total space to define connections
on a parabolic $G$--bundle; see \cite{Bi2} for the details.

Let
$$
0\, \longrightarrow\, \text{ad}(F_G)\otimes\Omega^1_Y
\, \longrightarrow\,\text{At}(F_G)\otimes\Omega^1_Y
\stackrel{q}{\longrightarrow}\,
TY\otimes\Omega^1_Y \longrightarrow\,0
$$
be the exact sequence obtained by tensoring \eqref{eq-1} with
$\Omega^1_Y$. Consider the inclusion of ${\mathcal O}_Y$
in $TY\bigotimes\Omega^1_Y$ obtained by sending the constant
function $1$ to the identity automorphism of $TY$.
Therefore, from the above exact sequence we get the following
short exact sequence of holomorphic vector bundles 
\begin{equation}\label{eq-2}
0\, \longrightarrow\, \text{ad}(F_G)\otimes\Omega^1_Y
\, \longrightarrow\,\widetilde{\text{At}}(F_G)\,:=\,
q^{-1}({\mathcal O}_Y) \, \stackrel{q}{\longrightarrow}\,
{\mathcal O}_Y \longrightarrow\,0
\end{equation}
over $Y$.

It is easy to see that giving a holomorphic (respectively
$C^\infty$) splitting of \eqref{eq-1} is equivalent to
giving a holomorphic (respectively
$C^\infty$) splitting of \eqref{eq-2}. Therefore, the twisted
version of the Atiyah exact sequence given in \eqref{eq-2}
is also suitable for defining connections.

The exact sequence in \eqref{eq-2} generalizes to the context
of parabolic $G$--bundles. Connections on a parabolic $G$--bundle
can be defined to be the splittings of the corresponding
short exact sequence.

\section{Preliminaries}\label{sec2}

Let $X$ be a connected smooth projective variety
of dimension $d$ defined over $\mathbb C$.
Let $D\, \subset\, X$ be a simple normal crossing hypersurface.
This means that $D$ is an effective and reduced
divisor with each irreducible component
of $D$ being smooth, and furthermore,
the irreducible components of $D$ intersect transversally. Let
\begin{equation}\label{decomp.-D}
D\, =\, \sum_{i=1}^\ell D_i
\end{equation}
be the decomposition of $D$ into irreducible components.
The above condition that the irreducible components of $D$
intersect transversally means that if
\begin{equation}\label{pt.-intsc.}
x\, \in\, D_{i_1}\cap D_{i_2}\cap \cdots \cap D_{i_k}\,
\subset\, D
\end{equation}
is a point where $k$ distinct components of $D$ meet, and $f_{i_j}$,
$j\,\in\, [1\, ,k]$, is the local equation of the divisor
$D_{i_j}$ around $x$, then $\{\text{d}f_{i_j}(x)\}$ is a linearly
independent subset of the holomorphic cotangent space $T^*_xX$ of
$X$ at $x$. This implies that for any choice of $k$ integers
\begin{equation}\label{k-int.}
1\, \leq \,  i_1 \, <\, i_2\, <\, \cdots\, < \, i_k \, \leq\,\ell\, ,
\end{equation}
each connected component of
$D_{i_1}\bigcap D_{i_2}\bigcap \cdots \bigcap D_{i_k}$ is a smooth
subvariety of $X$.

Let $E$ be an algebraic vector bundle over $X$. For
each $i\,\in\, [1\, , \ell]$, let
\begin{equation}\label{divisor-filt.}
E\vert_{D_i}\, =\, F^i_1 \,\supsetneq\, F^i_2 \,\supsetneq\,
F^i_3 \,\supsetneq\, \cdots \,\supsetneq\, F^i_{m_i} \,\supsetneq\,
F^i_{m_i+1}\, =\, 0
\end{equation}
be a filtration by subbundles of the restriction of $E$ to $D_i$.
In other words, each $F^i_j$ is a subbundle of $E\vert_{D_i}$
and $\text{rank}(F^i_j) \, >\, \text{rank}(F^i_{j+1})$
for $j\in [1\, , m_i]$.

A \textit{quasiparabolic} structure on $E$ over $D$ is a filtration
as above of each $E\vert_{D_i}$ satisfying the following extra
condition: Take any $k\,\in\, [1\, ,\ell]$,
and take integers $\{i_j\}_{j=1}^k$ as in \eqref{k-int.}.
If we fix some $F^{i_j}_{n_j}$, $n_j\, \in\,
[1\, , m_{i_j}]$, then over each connected component $S$
of $D_{i_1}\bigcap D_{i_2}\bigcap \cdots \bigcap D_{i_k}$,
the intersection
$$
\bigcap_{j=1}^k F^{i_j}_{n_j}\, \subset\,
E\vert_{D_{i_1}\cap \cdots \cap D_{i_k}}
$$
gives a subbundle of the restriction of $E$ to $S$.
It should be clarified that the rank of this
subbundle may depend on the choice of the component
$S\, \subset\, 
D_{i_1}\bigcap D_{i_2}\bigcap \cdots \bigcap D_{i_k}$.

To explain the above condition we give an example. Let $S\,
\subset\, D_{i_1}\bigcap D_{i_2}$ be a connected component,
where $1\, \leq \,  i_1 \, < \, i_2\, \leq\, \ell$, and take
$F^{i_1}_{n_1}$, $F^{i_2}_{n_2}$, where $n_j\, \in\,
[1\, , m_{i_j}]$, $j\,=\, 1,2$.
For any point $x\,\in\, S$, consider the subspace
$(F^{i_1}_{n_1})_x\bigcap (F^{i_2}_{n_2})_x\, \subset\, E_x$.
The above condition says that the dimension of this subspace is
independent of the choice of $x\in S$. But this dimension depends
on the choices of $i_j$, $n_j$, $j\, \in\, [1\, ,2]$, and it also
depends on the choice of the connected component $S$ in $D_{i_1}
\bigcap D_{i_2}$. Note that the condition that $(F^{i_1}_{n_1})_x
\bigcap (F^{i_2}_{n_2})_x$ is of constant dimension over a connected
component $S$ is equivalent to the condition that $(F^{i_1}_{n_1}
\bigcap F^{i_2}_{n_2})\vert_S$ is a subbundle of $E\vert_S$.

For a quasiparabolic structure as above, \textit{parabolic
weights} are a collection of rational numbers
\begin{equation}\label{divisor-para.-weight}
0\, \leq\, \lambda^i_1\, < \, \lambda^i_2\, < \,
\lambda^i_3 \, < \,\cdots \, < \, \lambda^i_{m_i}
\, <\, 1
\end{equation}
where $i\in [1\, ,\ell]$. The parabolic weight $\lambda^i_j$
corresponds to the subbundle $F^i_j$ in \eqref{divisor-filt.}. A
\textit{parabolic structure} on $E$ is a quasiparabolic
structure with parabolic weights. A vector bundle
over $X$ equipped
with a parabolic structure on it is also called a
\textit{parabolic vector bundle}.

For notational convenience, a parabolic vector bundle defined
as above will be denoted by $E_*$. The divisor $D$ is called
the \textit{parabolic divisor} for $E_*$.

Let $G$ be a connected complex linear algebraic group. We will
recall the definition of a parabolic $G$--bundle introduced
in \cite{BBN1}.

Let $\text{Rep}(G)$ denote the category of
all finite dimensional
rational left representations of $G$. Let $\text{Vect}(X)$
denote the category of algebraic vector bundles over $X$.
Nori showed that a principal $G$--bundle over $X$ is a
functor from $\text{Rep}(G)$ to $\text{Vect}(X)$ which is compatible
with the operations of taking direct sum, tensor product and dual.
Given a principal $G$--bundle $E_G$ over $X$, the corresponding
functor $\text{Rep}(G)\,\longrightarrow\, \text{Vect}(X)$
sends a $G$--module $V$ to the vector bundle
$E_G\times^GV$ over $X$ associated to $E_G$ for $V$;
see \cite{No1}, \cite{No2} for the details.

Let $\text{Pvect}(X)$ denote the category of
all parabolic vector bundles
over $X$ with $D$ as the parabolic divisor. In \cite{BBN1}, parabolic
$G$--bundles over $X$ with $D$ as the parabolic divisor were defined
to be as functors from $\text{Rep}(G)$ to $\text{Pvect}(X)$ satisfying
a list of conditions identical to the list of conditions of Nori
in the above mentioned characterization
of usual principal bundles as functors. It may be mentioned
the operations of taking direct sum, tensor product and
dual of usual vector bundles are replaced by the operations of
taking parabolic direct sum, parabolic tensor product and parabolic
dual. See \cite{BBN1} for the details.

In \cite{BBN2}, the notion of a ramified $G$--bundle over a curve
was introduced. The main result of \cite{BBN2} is the construction
of a natural bijective correspondence between the ramified
$G$--bundles over a Riemann surface $C$, with ramifications over a
finite set of points $D_0\, \subset\,
C$, and the parabolic $G$--bundles over $C$ with $D_0$
as the parabolic divisor.

We will define ramified $G$--bundles over the projective manifold $X$
with ramification over the simple normal crossing divisor $D$.

A \textit{ramified $G$--bundle} over $X$ with ramification over
$D$ is a smooth complex variety $E_G$ on which $G$ acts (algebraically)
on the right, that is, the map
$$
f \, :\, E_G\times G\, \longrightarrow\, E_G
$$
defining the action is an algebraic morphism, together with
a surjective algebraic map
\begin{equation}\label{def.-psi}
\psi\, :\, E_G\, \longrightarrow\, X
\end{equation}
satisfying the following five conditions:
\begin{enumerate}
\item{} $\psi\circ f \, =\, \psi\circ p_1$, where $p_1$ is the
natural
projection of $E_G\times G$ to $E_G$, that is;

\item{} for each point $x\, \in\, X$, the action of $G$ on the
reduced fiber $\psi^{-1}(x)_{\text{red}}$ is transitive;

\item{} the restriction of $\psi$ to $\psi^{-1}(X\setminus D)$ makes
$\psi^{-1}(X\setminus D)$ a principal $G$--bundle over
$X\setminus D$, that is, the map $\psi$ is smooth over
$\psi^{-1}(X\setminus D)$ and the map to the fiber product
$$
\psi^{-1}(X\setminus D)\times G\, \longrightarrow\,
\psi^{-1}(X\setminus D)\times_{X\setminus D} 
\psi^{-1}(X\setminus D)
$$
defined by $(z\, ,g) \longmapsto\, (z\, , f(z,g))$ is an
isomorphism;

\item{} for each irreducible component $D_i\, \subset\, D$,
the reduced inverse image $\psi^{-1}(D_i)_{\text{red}}$ is a
smooth divisor and
\begin{equation}\label{def.-wh.-D}
\widehat{D}\, :=\, \sum_{i=1}^\ell \psi^{-1}(D_i)_{\text{red}}
\end{equation}
is a normal crossing divisor on $E_G$;

\item{} for any smooth point $z\in \widehat{D}$, the isotropy
group $G_z\, \subset\, G$, for the action of $G$ on $E_G$,
is a finite cyclic group that acts faithfully on the
quotient line $T_zE_G/T_z\psi^{-1}(D)_{\text{red}}$.
\end{enumerate}

Note that since the map
$\psi$ commutes with the action of $G$, the
isotropy subgroup $G_z$ preserves $T_z\psi^{-1}(D)_{\text{red}}
\,\subset\, T_z E_G$.
Therefore, there is an induced action of $G_z$ on the fiber
$T_zE_G/T_z\psi^{-1}(D)_{\text{red}}$ of the normal bundle.
Since the finite isotropy
group $G_z$ in the above condition (5) acts faithfully
on the line $T_zE_G/T_z\psi^{-1}(D)_{\text{red}}$, it follows
automatically that $G_z$ is a cyclic group.

Let $E_G$ be a ramified $G$--bundle over $X$ with ramification
over the divisor $D$.
Fix a component $D_i\, \subset\, D$. Let $x\, \in\, D_i$ be
a smooth point of $D$. The order of the finite cyclic group
$G_z$, where $z$ satisfies the condition $\psi(z)\, =\, x$,
does not depend on the choices of $x$ and $z$; it depends
only on the component $D_i$ and $E_G$. Therefore, given any
ramified $G$--bundle $E_G$, we have a positive integer
$\eta_i$ associated to each component $D_i$; for any $z$
as above, the order of $G_z$ is $\eta_i$.

In the next section we will describe a correspondence
between parabolic $G$--bundles and ramified $G$--bundles. 

\section{Ramified $G$--bundles as parabolic
$G$--bundles}\label{sec3}

The map from ramified $G$--bundles to
parabolic $G$--bundles described in Section 2 of \cite{Bi2}.
Although it is assumed in \cite{Bi2} that $\dim_{\mathbb C}X
\,=\, 1$, the construction of a parabolic $G$--bundle from a
given ramified $G$--bundle goes through.

Let $E_*$ be a parabolic $G$--bundle over over $X$ with $D$ as
the parabolic divisor. In \cite{BBN1} the following was proved:

There is a Galois covering
\begin{equation}\label{varphi}
\varphi\, :\, Y\, \longrightarrow X
\end{equation}
and a $\Gamma$--linearized principal $G$--bundle $E_G$ over $Y$,
where $\Gamma$ is the Galois group for $\varphi$, such that
$E_*$ corresponds to $E_G$ (see \cite[Theorem 4.3]{BBN1}).

We will show that $E_G/\Gamma$ is a ramified $G$--bundle over $X$. All
the properties except one of a ramified $G$--bundle for $E_G/\Gamma$
has already been shown in \cite{BBN2}. The only property that remains
to be checked is that $E_G/\Gamma$ is smooth. (The argument given in
\cite{BBN2} that $E_G/\Gamma$ is smooth uses the assumption that
$\dim_{\mathbb C} X \, =\,1$.) We will show that $E_G/\Gamma$ is smooth.

Let $S\, \subset\, Y$ be the subscheme where the map $\varphi$ in
\eqref{varphi} fails to be smooth. Since $X\,=\, Y/\Gamma$ is smooth,
each component of $S$ is a hypersurface. For any rational point
$y\, \in\, S$, let $\Gamma_y\, \subset\, \Gamma$ be the isotropy
subgroup for the action of the Galois group $\Gamma$ on $Y$.

For any rational point $y\, \in\, S\setminus \varphi^{-1}(D)$ in the
complement of the inverse image $\varphi^{-1}(D)$,
the action of $\Gamma_y$ on the fiber $(E_G)_y$ is trivial. Indeed,
this follows from the fact the action of $G$ on the fiber
$(E_*)_{\varphi(y)}$ is a free action (recall the definition of
a ramified $G$--bundle with ramifications over $D$).

For each irreducible component $D_i$ of $D$, there is a maximal
subgroup $H_i\, \subset\, \Gamma$ such that $H_i\, \subset\,
\Gamma_y$ for all $y\, \in\, \varphi^{-1}(D_i)$. Furthermore,
for the general point $y\, \in\, \varphi^{-1}(D_i)$,
the equality $H_i\, =\, \Gamma_y$ holds, and also the group
$H_i$ is cyclic. The action of $H_i$ on the fiber of $E_G$
over a point $\varphi^{-1}(D_i)$ need not be free,
but there is a fixed quotient group $H'_i$ of $H_i$
such that for all $y\, \in\, \varphi^{-1}(D_i)$, the
action of $H_y$ on the fiber $(E_G)_y$ factors through
$H'_i$. Furthermore, the action of $H'_i$ on $(E_G)_y$
is free. From these it follows that the quotient
$E_G/\Gamma$ is smooth.

\section{Holomorphic connections}

We will first recall the usual Atiyah exact sequence and its
equivalent formulations.

\subsection{The Atiyah exact sequence}

Let $M$ be a complex manifold and
\begin{equation}\label{de2f}
f\, :\, E_H\, \longrightarrow\, M
\end{equation}
a holomorphic principal
$H$--bundle over $M$, where $H$ is a complex Lie group. Consider
the sheaf $\mathcal F$
on $M$ that associates to any open subset $U\, \subset\,
M$ the space of all $H$--invariant holomorphic vector fields on
$f^{-1}(U)\, \subset\, E_H$. Therefore, ${\mathcal F}(U)$ is a
${\mathcal O}_U$--module, where ${\mathcal O}_U$ is the algebra
of holomorphic functions on $U$; the multiplication of a
holomorphic vector field $\tau$ on $f^{-1}(U)$ with a function
$\phi\, \in\, {\mathcal O}_U$ is the vector field
$(\phi\circ f)\cdot  \tau$. Since the action of $H$ on the fibers
of $f$ is transitive, it follows immediately that ${\mathcal F}$
is a locally free coherent analytic sheaf on $M$.

The holomorphic vector bundle over $M$ defined by $\mathcal F$ is
called the \textit{Atiyah bundle} for $E_H$. The Atiyah bundle for
$E_H$ is denoted by $\text{At}(E_H)$.

Let $\text{ad}(E_H)$ be the adjoint bundle of $E_H$. So
$\text{ad}(E_H)$ is the holomorphic vector bundle over $M$
associated to $E_H$ for the adjoint action of $H$ on its Lie
algebra $\mathfrak h$. We recall that $\mathfrak h$ is identified
with the vector fields on $H$ invariant under the right translation
action of $H$ on itself. Using this it follows that for any
open subset $U\, \subset\, M$, the space of all holomorphic sections
of $\text{ad}(E_H)$ over $U$ is identified with the space of all
$H$--invariant holomorphic vector fields on
$f^{-1}(U)\, \subset\, E_H$ that lie in the kernel of the differential
$$
{\rm d}f \, :\, TE_H\, \longrightarrow\, f^*TM
$$
of the projection $f$ in \eqref{de2f}. In other words,
$\Gamma (U,\, \text{ad}(E_H))$ is the space of all $H$--invariant
holomorphic vertical vector fields on $f^{-1}(U)$.
Therefore, we obtain
a short exact sequence of holomorphic vector bundles
\begin{equation}\label{at.ex.}
0\, \longrightarrow\, \text{ad}(E_H)\, \longrightarrow\,
\text{At}(E_H) \longrightarrow\, TM \longrightarrow\,0
\end{equation}
over $M$; the projection $\text{At}(E_H) \longrightarrow\, TM$
is give by the differential ${\rm d}f$.

The following definitions are from \cite{At}.

\begin{definition}\label{def.c0}
{\rm A \textit{holomorphic connection} on $E_H$ is a holomorphic
splitting of the short exact sequence in \eqref{at.ex.}. A
\textit{complex connection} on $E_H$ is a $C^\infty$ splitting
of the short exact sequence in \eqref{at.ex.}.}
\end{definition}

Let
\begin{equation}\label{at0.ex.}
0\, \longrightarrow\, \text{ad}(E_H)\otimes
\Omega^1_M\, \longrightarrow\, \text{At}(E_H)\otimes \Omega^1_M \,
\stackrel{q}{\longrightarrow}\,
TM\otimes \Omega^1_M \, \longrightarrow\,0
\end{equation}
be the short exact sequence obtained by tensoring \eqref{at.ex.}
with the holomorphic cotangent bundle $\Omega^1_M$. Using the section
of $TM\otimes T^*M$ given by the identity automorphism of $TM$, the
structure sheaf ${\mathcal O}_M$ is a subsheaf of $TM\otimes
\Omega^1_M$.
Therefore, from \eqref{at0.ex.} we have the short exact sequence of
holomorphic vector bundles
\begin{equation}\label{at.ex.2}
0\, \longrightarrow\, \text{ad}(E_H)\otimes\Omega^1_M\,\longrightarrow
\, \widetilde{\text{At}}(E_H)\, :=\, q^{-1}({\mathcal O}_M)
\, \stackrel{q\vert_{\widetilde{\text{At}}(E_H)}}{\longrightarrow}\,
{\mathcal O}_M \, \longrightarrow\,0
\end{equation}
over $M$, where $q$ is the projection in \eqref{at0.ex.}.

\begin{remark}\label{c.u}
{\rm A holomorphic splitting of \eqref{at.ex.} gives a holomorphic
splitting of \eqref{at.ex.2} and conversely a holomorphic splitting
of \eqref{at.ex.2} gives a holomorphic splitting of \eqref{at.ex.}.
Similarly, giving a $C^\infty$
splitting of the exact sequence in \eqref{at.ex.} is equivalent to
giving a $C^\infty$ splitting of \eqref{at.ex.2}. Therefore, giving
a holomorphic connection on $E_G$ is equivalent to giving a
holomorphic splitting of \eqref{at.ex.2}. Similarly, giving
a complex connection on $E_G$ is equivalent to giving a
$C^\infty$ splitting of \eqref{at.ex.2}.}
\end{remark}

\subsection{Connections on a ramified principal bundle}

Let $E_G$ be a ramified $G$-bundle over $X$ with ramification
over $D$. As in \eqref{def.-psi}, let $\psi$ denote the projection
of $E_G$ to $X$. Consider the subbundle of the holomorphic
tangent bundle
$$
{\mathcal K} \, \subset\, TE_G
$$
defined by the orbits of the action of
$G$ on $E_G$. Since the isotropy
subgroups, for the action of $G$ on $E_G$,
are all finite subgroups of $G$, it follows that ${\mathcal K}$
is a subbundle of $TE_G$, and the vector bundle ${\mathcal K}$ is
identified with the trivial vector bundle over $E_G$ with fiber
$\mathfrak g$, where $\mathfrak g$ is the Lie algebra of $G$.
The differential
$$
{\rm d}\psi\, :\, TE_G\, \longrightarrow\, \psi^*TX
$$
evidently vanishes on $\mathcal K$.

Let ${\mathcal Q}$ denote the quotient bundle $TE_G/{\mathcal K}$.
So we have a short exact sequence of holomorphic vector bundles
\begin{equation}\label{ex.p.0}
0\, \longrightarrow\, {\mathcal K}\, \longrightarrow\,
TE_G\, \longrightarrow\, {\mathcal Q} \, \longrightarrow\,0
\end{equation}
over $E_G$. Tensoring \eqref{ex.p.0} with ${\mathcal Q}^*$
we get the exact sequence
\begin{equation}\label{ex.p.1}
0\, \longrightarrow\, {\mathcal K}\otimes {\mathcal Q}^*\,
\longrightarrow\,
TE_G\otimes {\mathcal Q}^*\, \stackrel{q_0}{\longrightarrow}\,
{\mathcal Q}\otimes{\mathcal Q}^* \, \longrightarrow\,0
\end{equation}
over $E_G$. As in \eqref{at.ex.2}, we will consider the
inverse image of the trivial line sub--bundle of
${\mathcal Q}\otimes{\mathcal Q}^*$ generated by the identity
automorphism of ${\mathcal Q}$. So we have the
short exact sequence of holomorphic vector bundles
\begin{equation}\label{ex.p}
0\, \longrightarrow\, {\mathcal K}\otimes {\mathcal Q}^*\,
\longrightarrow\, {\mathcal V}_{E_G}\, :=\,
q^{-1}_0({\mathcal O}_{E_G})\,
\stackrel{q_0\vert_{{\mathcal V}_{E_G}}}{\longrightarrow}
\, {\mathcal O}_{E_G} \, \longrightarrow\, 0
\end{equation}
over $E_G$.

We note that the action of $G$ on $E_G$ has natural lifts to all
the three vector bundles in the short exact sequence in \eqref{ex.p}.
Furthermore, all the homomorphisms in \eqref{ex.p} commute with
the actions of $G$. Therefore, the direct image, on $X$, of any of
the vector bundles in \eqref{ex.p} is equipped with an action of $G$.

Define the quasi--coherent analytic sheaves
\begin{equation}\label{dag}
{\mathcal A}_{E_G}\, :=\,
(\psi_*({\mathcal K}\otimes {\mathcal Q}^*))^G
\end{equation}
and
\begin{equation}\label{dbg}
{\mathcal B}_{E_G} \, :=\, (\psi_*{\mathcal V}_{E_G})^G
\end{equation}
on $X$, where $\psi$ is the projection in \eqref{def.-psi},
and by $W^G$, where $W$ is any sheaf on $X$ equipped
with an action of $G$, we mean the $G$--invariant
part of $W$. Since the action of $G$ on the
fibers of $\psi$ is transitive with finite isotropy
subgroups, it follows that
both ${\mathcal A}_{E_G}$ and ${\mathcal B}_{E_G}$
introduced in \eqref{dag}
and \eqref{dbg} are locally free coherent analytic sheaves on $X$.
The holomorphic vector bundles over $X$ defined by ${\mathcal A}_{E_G}$
and ${\mathcal B}_{E_G}$ will also be denoted by
${\mathcal A}_{E_G}$ and ${\mathcal B}_{E_G}$ respectively.

We note that the action of $G$ on the sheaf ${\mathcal O}_{E_G}$ 
in \eqref{ex.p} is the trivial action. This means that the identity
automorphism of $\mathcal Q$ is preserved by the action of $G$ on
${\mathcal Q}\otimes{\mathcal Q}^*$. Therefore,
$(\psi_* {\mathcal O}_{E_G})^G\, =\, {\mathcal O}_X$.

Using the above observations, from \eqref{ex.p} we have following
short exact sequence of holomorphic vector bundles over $X$
\begin{equation}\label{ar}
0\, \longrightarrow\, {\mathcal A}_{E_G}\,
\longrightarrow\, {\mathcal B}_{E_G}\, \longrightarrow\,
{\mathcal O}_X \, \longrightarrow\,0\, .
\end{equation}

\begin{definition}\label{def2}
{\rm A \textit{holomorphic connection} on $E_G$ is a holomorphic
splitting of the short
exact sequence in \eqref{ar}. A \textit{complex
connection} on $E_G$ is a $C^\infty$
splitting of the short exact sequence in \eqref{ar}.}
\end{definition}

When $E_G$ is an usual principal $G$--bundle, the exact sequence in
\eqref{ar} clearly coincides with the exact sequence in \eqref{at.ex.2}.
In view of Remark \ref{c.u}, the above definitions coincide with
those in Definition \ref{def.c0} when $E_G$ is an usual principal
bundle.

In \cite{Bi2} we defined connections on a ramified $G$--bundle
over a curve. The following theorem shows that the
above definition coincides with the one given in \cite{Bi2}.

\begin{theorem}\label{thm1}
Let $E_G$ be a ramified
principal $G$--bundle over $X$. Giving a holomorphic connection
on $E_G$ is equivalent to giving a holomorphic connection on $E_G$
in the sense of \cite{Bi2}. Similarly, giving a complex connection
on $E_G$ is equivalent to giving a complex connection on $E_G$
in the sense of \cite{Bi2}.
\end{theorem}

\begin{proof}
Let $\beta\, :\, {\mathcal O}_X \, \longrightarrow\,
{\mathcal B}_{E_G}$ be a holomorphic splitting of the exact sequence
in \eqref{ar}. The lift on $\beta$ to $E_G$ gives a
$G$--equivariant holomorphic splitting
\begin{equation}\label{wtb}
\widetilde{\beta}\, :\, {\mathcal O}_{E_G}\, \longrightarrow\,
{\mathcal V}_{E_G}
\end{equation}
of the exact sequence in \eqref{ex.p}. Therefore,
$\widetilde{\beta}$ gives a homomorphism of holomorphic vector
bundles
$$
{\beta}'\, :\,  {\mathcal Q}\, \longrightarrow\, TE_G
$$
whose composition with the projection
$TE_G\, \longrightarrow\, {\mathcal Q}$ in \eqref{ex.p.0}
is the identity automorphism of ${\mathcal Q}$. Let
\begin{equation}\label{wtga}
\widetilde{\gamma}\, :\, TE_G\, \longrightarrow\,
{\mathcal K}
\end{equation}
be the projection given by the above homomorphism ${\beta}'$,
where $\mathcal K$ is the kernel in \eqref{ex.p.0}. Therefore,
the kernel of $\widetilde{\gamma}$ is the image of $\beta'$,
and the composition of $\widetilde{\gamma}$ with the inclusion
${\mathcal K}\, \hookrightarrow\, TE_G$ in \eqref{ex.p.0}
is the identity automorphism of ${\mathcal K}$.

We recall that ${\mathcal K}$ is the trivial vector bundle over
$E_G$ whose fiber is the Lie algebra $\mathfrak g$ of $G$. Therefore,
the homomorphism $\widetilde{\gamma}$ in \eqref{wtga} defines a
$\mathfrak g$--valued holomorphic one--form $\gamma$ on $E_G$.
Since the homomorphism $\widetilde{\beta}$ in \eqref{wtb} is
$G$--equivariant, we conclude that the holomorphic one--form
$\gamma$ on $E_G$ is also $G$--equivariant for the adjoint action
of $G$ on $\mathfrak g$. From the fact that the composition of
$\widetilde{\gamma}$ with the inclusion
${\mathcal K}\, \hookrightarrow\, TE_G$ in \eqref{ex.p.0} is
the identity automorphism of ${\mathcal K}$ it follows immediately
that the restriction of $\gamma$ to each fiber of the projection
$\psi$ (see \eqref{def.-psi}) is the Maurer--Cartan form on the
fiber. Therefore, $\gamma$ defines a holomorphic connection on
$E_G$ in the sense of \cite{Bi2}.

Conversely, let $\gamma$ be a $\mathfrak g$--valued holomorphic
one--form on $E_G$ defining a holomorphic connection on $E_G$
in the sense of \cite{Bi2}. Therefore, $\gamma$ is $G$--equivariant
and it coincides with the Maurer--Cartan form on the fibers
of the projection $\psi$. Consequently, $\gamma$ gives a
$G$--invariant holomorphic splitting
$$
\widetilde{\gamma}\, :\, TE_G\, \longrightarrow\,
{\mathcal K}
$$
of the exact sequence in \eqref{ex.p.0}. Consider the corresponding
homomorphism
$$
{\beta}'\, :\,  {\mathcal Q}\, \longrightarrow\, TE_G\, .
$$
So, the image of ${\beta}'$ is the kernel of $\widetilde{\gamma}$,
and the composition of ${\beta}'$ with the projection
$TE_G\, \longrightarrow\, {\mathcal Q}$ in \eqref{ex.p.0}
is the identity automorphism of ${\mathcal Q}$. Therefore,
$\beta'$ gives a $G$--invariant holomorphic splitting
$$
\widetilde{\beta}\, :\, {\mathcal O}_{E_G}\, \longrightarrow\,
{\mathcal V}_{E_G}
$$
of the exact sequence in \eqref{ex.p}. Hence
$\widetilde{\beta}$ descends to a holomorphic splitting of the
exact sequence in \eqref{ar}.

Therefore, giving a holomorphic connection
on $E_G$ is equivalent to giving a holomorphic connection on $E_G$
in the sense of \cite{Bi2}.

Similarly, it can be shown that giving a complex connection on
$E_G$ is equivalent to giving a complex connection on $E_G$
in the sense of \cite{Bi2}. This completes the proof of the theorem.
\end{proof}

\subsection{A construction of ${\mathcal A}_{E_G}$}

Consider the derivation action of $TX$ on ${\mathcal O}_X$.
Let
$$
TX(-\log D)\, \subset\, TX
$$
be the subsheaf that leaves ${\mathcal O}_X(-D)\, \subset\,
{\mathcal O}_X$ invariant. So
$$
TX\otimes {\mathcal O}_X(-D)\, \subset\, TX(-\log D)
\, \subset\, TX\, ,
$$
and if $x\, \in\, D_i$ is a smooth point of $D$, then the
image of the fiber $TX(-\log D)_x$ in $T_xX$ coincides 
with $T_xD_i$. More generally, for any
point $x$ as in \eqref{pt.-intsc.}, the
image of the fiber $TX(-\log D)_x$ in $T_xX$ is contained
in the kernel of the projection of $T_xX$ to the
fiber over $x$
of the normal bundle to $D_{i_1}\bigcap D_{i_2}
\bigcap\cdots \bigcap D_{i_k}$. The subsheaf $TX(-\log D)$ is
locally free and hence it defines a holomorphic
vector bundle over $X$. The dual vector bundle
$(TX(-\log D))^*$ is denoted by $\Omega^1_X(\log D)$.

Let $E_G$ be a ramified $G$--bundle over $X$ with
ramification over $D$. As we saw in Section \ref{sec3}, the
ramified $G$--bundle $E_G$ gives a functor from $\text{Rep}(G)$
to $\text{PVect}(X)$. Consider the $G$--module $\mathfrak g$
equipped with the adjoint action of $G$. The image of
this $G$--module $\mathfrak g$ by the functor
$\text{Rep}(G)\,\longrightarrow\, \text{PVect}(X)$ corresponding
to $E_G$ will be denoted by $E^{\mathfrak g}_*$. So
$E^{\mathfrak g}_*$ is a parabolic vector bundle over
$X$ with parabolic structure over $D$.

Let $E^{\mathfrak g}_0$ denote the underlying holomorphic
vector bundle for the above defined parabolic vector bundle
$E^{\mathfrak g}_*$. For each irreducible component $D_i$
of $D$, let
$$
E^{\mathfrak g}_0\vert_{D_i}\, =\, F^i_1 \,\supset\, F^i_2 \,
\supset\, F^i_3 \,\supset\, \cdots \,\supset\, F^i_{m_i} \,
\supset\, F^i_{m_i+1}\, =\, 0
$$
be the filtration as in \eqref{divisor-filt.}. We will define a
subbundle $V^i$ of $E^{\mathfrak g}_0\vert_{D_i}$.
If the parabolic
weight $\lambda_1$ corresponding to $F^i_1$ is zero, then set
$V^i\, :=\, F^i_2$, and if $\lambda_1\, \not=\, 0$, then set
$V^i\, :=\, F^i_1$. Let ${\mathcal F}\, \subset\,
E^{\mathfrak g}_0$ be the subsheaf that fits in the following
short exact sequence of coherent sheaves on $X$
\begin{equation}\label{st.ex.-1}
0\, \longrightarrow\, {\mathcal F}\, \longrightarrow\,
E^{\mathfrak g}_0\, \longrightarrow\,
\bigoplus_{i=1}^\ell (E^{\mathfrak g}_0\vert_{D_i})/V_i
\, \longrightarrow\,0\, ,
\end{equation}
where $V_i$ are defined above, and $\{D_i\}_{i=1}^\ell$ are
the irreducible components of $D$. Therefore, ${\mathcal F}$
is a holomorphic vector bundle over $X$ which is identified with
$E^{\mathfrak g}_0$ over the complement $X\setminus D$.

Let
\begin{equation}\label{st.ex.0}
{\mathcal W}_{E_G}\, \subset\, 
E^{\mathfrak g}_0\otimes \Omega^1_X(\log D)
\end{equation}
be the coherent subsheaf generated by
the two subsheaves ${\mathcal F}\bigotimes
\Omega^1_X(\log D)$ and $E^{\mathfrak g}_0\bigotimes\Omega^1_X$,
where $\mathcal F$ is defined in \eqref{st.ex.-1} and
$\Omega^1_X(\log D)$ was defined earlier. It is easy to check that
${\mathcal W}_{E_G}$ is locally free. Hence it defines a holomorphic
vector bundle over $X$. This vector bundle is clearly
identified with $E^{\mathfrak g}_0\bigotimes\Omega^1_X$
over the complement $X\setminus D$.

\begin{proposition}\label{prop1}
The vector bundle ${\mathcal A}_{E_G}$ in \eqref{ar} is identified
with the vector bundle ${\mathcal W}_{E_G}$ in \eqref{st.ex.0}.
\end{proposition}

\begin{proof}
There is a Galois covering
\begin{equation}\label{f.l}
f\, : Y\, \longrightarrow\, X
\end{equation}
with $Y$ a smooth complex projective variety and
a holomorphic principal $G$--bundle $F_G$ over $Y$ equipped
with a lift of the action of the Galois group $\Gamma\, :=\,
\text{Gal}(f)$ on $Y$ such that $E_G \, =\, F_G/\Gamma$. Indeed,
in \cite{BBN1} it was shown that given a parabolic $G$--bundle
over $X$, such a pair $(f\, , F_G)$ exist. We noted in Section
\ref{sec3} that a parabolic $G$--bundle is same as a
ramified $G$--bundle.

Using the actions of the Galois group $\Gamma$ on $Y$ and $F_G$,
the vector bundle $\text{ad}(F_G)\bigotimes \Omega^1_Y$ over $Y$
is equipped with an action of $\Gamma$. The corresponding parabolic
vector bundle over $X$ has the property that its underlying
vector bundle is identified with the vector bundle
${\mathcal W}_{E_G}$ defined in \eqref{st.ex.0}. (See \cite{Bi1}
for the correspondence between parabolic vector bundles over $X$
and $\Gamma$--linearized vector bundles over $Y$.)

Let $\psi_0\, :\, F_G\, \longrightarrow\, Y$ be the
natural projection. Let
$$
\delta\, :\, F_G\, \longrightarrow\, F_G/\Gamma\, =\, E_G
$$
be the quotient map.
Since the actions of $G$ and $\Gamma$ on $F_G$ commute, for any
$\Gamma$--linearized vector bundle $W$ over $Y$, the pullback
$\psi^*_0 W$ on $F_G$ has the following property: the invariant
direct image $(f_*W)^\Gamma$ on $X$, where $f$ is the projection
in \eqref{f.l}, is identified with the invariant direct
image $(\psi_*(\delta_*\psi^*_0 W)^\Gamma)^G$, where $\psi$
is the projection in \eqref{def.-psi}. Now setting $W$ to be
$\text{ad}(F_G)\bigotimes \Omega^1_Y$ we conclude that
the vector bundle
${\mathcal W}_{E_G}$ defined in \eqref{st.ex.0} is identified
with the vector bundle ${\mathcal A}_{E_G}$ in \eqref{ar}.
This completes the proof of the proposition.
\end{proof}


\end{document}